\documentclass[12pt,usletter]{article}
\usepackage[textwidth=6.5in,textheight=21.5cm]{geometry}
\usepackage{natbib, graphicx, setspace, amssymb, amsmath, amsthm, amstext, booktabs, fancyhdr, multirow, float, bm, xr, mathptmx, subfigure, microtype, url, dsfont, paralist}
\usepackage[dvipsnames]{xcolor}
\usepackage{placeins} 
\usepackage{authblk} 
\graphicspath{{plots/}}
\allowdisplaybreaks 

\newcommand*\patchAmsMathEnvironmentForLineno[1]{%
  \expandafter\let\csname old#1\expandafter\endcsname\csname #1\endcsname
  \expandafter\let\csname oldend#1\expandafter\endcsname\csname end#1\endcsname
  \renewenvironment{#1}%
     {\linenomath\csname old#1\endcsname}%
     {\csname oldend#1\endcsname\endlinenomath}}%
\makeatletter
\g@addto@macro\normalsize{%
  \setlength\abovedisplayskip{8pt}
  \setlength\belowdisplayskip{8pt}
  \setlength\abovedisplayshortskip{8pt}
  \setlength\belowdisplayshortskip{8pt}
}
\makeatother

\theoremstyle{plain}

 \newtheorem{theorem}{Theorem}  
\newtheorem{lemma}{Lemma}

\begin{document}
\title{Semiparametric Regression using Variational Approximations}
\author[1]{Francis K. C. Hui\thanks{Corresponding author: Mathematical Sciences Institute, The Australian National University, 2601, Acton, ACT, Australia. P: +612-6125-0581 E: francis.hui@anu.edu.au}} 
\author[2]{C. You}
\author[3]{H. L. Shang}
\author[4]{Samuel M\"{u}ller} 
\affil[1]{Mathematical Sciences Institute, The Australian National University, Canberra, Australia}
\affil[2]{School of Mathematical Sciences, University of Nottingham, Ningbo, China}
\affil[3]{Research School of Finance, Actuarial Studies, and Statistics, The Australian National University, Canberra, Australia}
\affil[4]{School of Mathematics and Statistics, University of Sydney, Sydney, Australia}
\date{}
\maketitle
\thispagestyle{empty} 

\begin{abstract}
Semiparametric regression offers a flexible framework for modeling non-linear relationships between a response and covariates. A prime example are generalized additive models where splines (say) are used to approximate non-linear functional components in conjunction with a quadratic penalty to control for overfitting. Estimation and inference are then generally performed based on the penalized likelihood, or under a mixed model framework. The penalized likelihood framework is fast but potentially unstable, and choosing the smoothing parameters needs to be done externally using cross-validation, for instance. The mixed model framework tends to be more stable and offers a natural way for choosing the smoothing parameters, but for non-normal responses involves an intractable integral. 

In this article, we introduce a new framework for semiparametric regression based on variational approximations. The approach possesses the stability and natural inference tools of the mixed model framework, while achieving computation times comparable to using penalized likelihood. Focusing on generalized additive models, we derive fully tractable variational likelihoods for some common response types. We present several features of the variational approximation framework for inference, including a variational information matrix for inference on parametric components, and a closed-form update for estimating the smoothing parameter. We demonstrate the consistency of the variational approximation estimates, and an asymptotic normality result for the parametric component of the model. Simulation studies show the variational approximation framework performs similarly to and sometimes better than currently available software for fitting generalized additive models. 

\end{abstract}

\section{Introduction} \label{sec:intro}

Semiparametric regression refers to a broad class of models which extend parametric regression so that part of the mean response is modeled as smooth functions of the covariates \citep[][]{ruppert2003semiparametric}. Because of their capacity to model complicated, non-linear relationships in a data-driven way, semiparametric regression is increasingly used in many fields of applied statistics. One well-known class of semiparametric regression models is generalized additive models \citep[GAMs,][]{hastie86}, where the linear predictor comprises a sum of smoothing functions with each smoothing function constructed using, for instance, penalized splines \citep{eilers1996}. Penalized splines combine a set of spline basis functions with a quadratic penalty on the corresponding smoothing coefficients in order to strike a balance between overfitting and accuracy. To estimate the GAM, a common approach is to then formulate it as a generalized linear model (GLM) augmented with the quadratic penalty, and then maximize the subsequent penalized likelihood function \citep{wood2006generalized}. This method is analogous to penalized quasi-likelihood (PQL) estimation \citep{breslow1993approximate}, which is often used to fit generalized linear mixed models. As with PQL estimation however, the penalized likelihood formulation can produce biased parameter estimates, particularly for discrete responses, and it is also known to suffer from potential convergence issues \citep{lin1999inference, wood2011fast}.

An alternative approach to estimation and inference for GAMs is to regard the penalty as a random effects distribution on the smoothing coefficients, thus allowing us to utilize the wealth of methods developed for mixed models and apply them to semiparametric regression. For instance, we can maximize the resulting marginal likelihood (defined as the likelihood of the responses integrated with respect to the random effects distribution; see Section~\ref{sec:gams} for an explicit mathematical definition) or an approximation thereof, such as the Laplace approximated likelihood \citep{kauermann2009some}. For inference, we can then use the posterior distribution of the smoothing coefficients as the basis for constructing hypothesis tests and prediction intervals \citep{wood2006generalized,krivobokova2010}, and use the marginal log-likelihood for smoothing parameter selection \citep{kauermann2005note,wood2016}. Empirical studies have also shown that the mixed model formulation for GAMs tends to be numerically more stable and produce more accurate estimates than the penalized likelihood formulation \citep{reiss2009smoothing,wood2011fast}. However, one major drawback of the mixed model approach is that the marginal likelihood does not have a tractable form for non-normal responses. This makes estimation and inference challenging, particularly when there are a large number of covariates and/or each covariate is modeled by a large number of basis functions. Indeed, the intractability of the integral is largely why the penalized likelihood formulation to estimation and inference remains the predominant approach for fitting semiparametric regression models e.g., both the \texttt{R} packages \texttt{mgcv} \citep{wood2006generalized} and \texttt{gamlss} \citep{rigby05} fit GAMs using the penalized likelihood formulation.

In this article, we introduce a new framework for estimation and inference for semiparametric regression based on variational approximations \citep[VA,][]{ormerod2010explaining}. Briefly, VA refers to a class of procedures aimed at obtaining a fully (or at least close to fully) tractable approximation to an intractable loss function. In the case of GAMs, we propose applying VA to obtain a fully closed-form lower bound to the mixed model marginal likelihood. It is important to highlight here the difference between VA and the increasingly popular variational Bayes (VB) method. Like standard Bayesian inference, VB assumes prior distributions for all parameters in the model (including the smoothing parameters) and aims to approximate the joint posterior distribution by a closed form, usually under simplifying assumptions such as mean field theory. For semiparametric regression, VB has been considered by a number of authors \citep[e.g.,][]{luts2014real, luts2015variational}, most of which have focused on developing efficient algorithms for estimation with little being explored in terms of theory and inference. By contrast, like the penalized likelihood and mixed model approaches reviewed above, our usage of the term ``variational approximation'', which is borrowed directly from \citet{hall2011asymptotic} and \citet{ormerod2012gaussian} among others, places it within the likelihood-based framework. That is, we use VA to refer to replacing the intractable marginal log-likelihood by a tractable lower bound approximation, which is then treated as the new objective function. The VA approach has only recently been proposed for overcoming intractable marginal likelihood functions, e.g., for generalized linear mixed models \citep{ormerod2012gaussian} and generalized linear latent variable models \citep{hui2016variational}. However to our knowledge, this article is the first to apply VA to semiparametric regression models, let alone study their theoretical properties and finite sample performance. As we will show, the VA approach to GAMs offers a balance between the stability and natural inference tools of the mixed model framework, while achieving computation times comparable to the penalized likelihood framework.

We apply VA to derive fully tractable variational likelihoods for the common cases of Poisson, normal, and Bernoulli responses, and show how these facilitate straightforward estimation of GAMs. Furthermore, we demonstrate how the approach naturally produces key inference tools such as confidence/prediction intervals and hypothesis tests for the smoothing coefficients. For the parametric component, we show how inference can be conducted by utilizing Louis' method \citep{louis82} to construct a fast, variational version of the observed Fisher information matrix. Regarding the critical choice of smoothing parameters, we show how the VA framework produces fully closed-form updates that are simple to compute.

Under mild regularity conditions, and assuming the true model is a GAM with a growing number of spline basis functions for each covariate, we establish consistency of the VA estimates and asymptotic normality for the parametric component of the model. These large-sample results are an important contribution on two fronts. First, relatively little has been done on the asymptotic behavior of VA, let alone of VB: \citet{hall2011asymptotic} showed consistency and asymptotic normality for the specific case of a Poisson mixed model with a random intercept, while \citet{ormerod2012gaussian} extended this to cover generalized linear mixed models involving a random intercept. In contrast, the results we develop in this article permit a growing number of spline basis functions i.e., multivariate random effects distributions. Second, much of the asymptotic theory developed for semiparametric regression has either assumed the underlying model can only be approximated by spline functions \citep[e.g.,][]{kauermann2009some,wang2011asymptotics}, or has assumed the true model consists of a fixed number of spline basis functions \citep[e.g.,][]{yu2002penalized,lu2008penalized}. The theory we develop here is the first to consider a setting which can be regarded as positioned between these two: we assume the true model consists of spline basis functions, but allow the number of basis functions for each covariate to grow polynomially with the sample size to allow increasing flexibility in modeling the response-covariate relationships. For establishing the asymptotic theory for GAMs, we believe that such a setting is very appropriate: research has shown that the number of basis functions needed to accurately capture the true mean response needs to only grow at a polynomial rate compared to sample size. That is, provided the number of knots is sufficiently large, the bias from approximating an underlying smooth function by a spline is negligible compared to both the bias due to a non-zero smoothing parameter and variance of the fitted curve \citep{ruppert2002,wang2011asymptotics}.

Simulation studies show that VA is substantially faster than methods which treat GAMs as mixed models and maximize the marginal likelihood, while also being more stable in small samples compared to methods which estimate GAMs by penalized likelihood. We apply the VA approach for fitting GAMs to uncover important associations between union membership and worker's attributes, with results showing evidence of a non-linear relationship between probability of union membership and hourly wages. 

To summarize, the main contributions of this article are: 
\begin{inparaenum}
\item[1)] Focusing on GAMs, we propose a new VA framework for estimation and inference in semiparametric regression; 
\item[2)] We demonstrate how key aspects of inference in GAMs can be performed naturally under this framework, including a new variational observed information matrix for parametric inference and fully closed-form updates for the smoothing parameters; 
\item[3)] We show that the estimates obtained via VA for a GAM are consistent and asymptotically normally distributed; 
\item[4)] Simulations demonstrate the strong performance, computational speed, and stability of the VA framework.
\end{inparaenum}

\section{Generalized Additive Models} \label{sec:gams}

For a set of independent observations $i = 1,\ldots,n$, let $y_i$ denote a univariate response, and suppose we have a set of predictors which can be split into a $p$-vector $\bm{x}_i = (x_{i1},\ldots,x_{ip})^{\top}$ of covariates that will be modeled parametrically, and a $q$-vector $\bm{u}_i = (u_{i1},\ldots,u_{iq})^{\top}$ of covariates that will be modeled nonparametrically. Conditional on all $p + q$ predictors, the $y_i$'s are assumed to be independent observations from the exponential family of distributions: $f(y_i | \cdot) = \exp[\phi^{-1}\{y_i \omega_i - b(\omega_i)\} + c(y_i,\phi)]$, where $b(\cdot)$ and $c(\cdot,\cdot)$ are known functions, $\omega_i$ is the canonical parameter, and $\phi$ is the dispersion parameter. For GAMs, the mean of the distribution $\mu_i = b'(\omega_i)$ is modeled as $g(\mu_i) = \eta_i = \bm{x}^{\top}_i \bm{\kappa} + \sum_{j=1}^q s_j(u_{ij})$, for some known link function $g(\cdot)$, where $\bm{\kappa}$ are the coefficients corresponding to the parametric component of the model (including the intercept) and $\{s_j(\cdot); j = 1,\ldots,q\}$ is a set of univariate smoothers for the semiparametric component. The above formulation is the classical form of a GAM involving additivity of the univariate smoothing functions for each covariate. Since the main contribution of this work is on an estimation and inference framework applicable to many types of semiparametric regression, then for elucidative reasons and as a first development into this field, we choose to work with this formulation of GAMs defined above.

One popular way of modeling the smoothing functions $s_j(\cdot)$ is through penalized splines, and here we focus on penalized B-splines or P-splines \citep{eilers1996}. For covariate $j = 1,\ldots,q$, suppose we choose a B-spline of degree $m_j$ and consider $K_j$ interior knots equally spaced on an interval which, without loss of generality, is taken to be $[0,1]$. Then we can write $s_j(u_{ij}) = \bm{z}^{\top}_{ij}\bm{\beta}_j$ where $\bm{z}_{ij}$ is constructed based on a set of $d_j = K_j+m_j-1$ basis functions evaluated at $u_{ij}$, and $\bm{\beta}_j$ is the corresponding $d_j$-vector of smoothing coefficients. Note the ``$-1$'' in the formula for $d_j$, which results from a centering constraint, which ensures that each smooth is centered around zero and the parameters are identifiable; details and references on the construction of P-splines are given in the Supplementary Material. To avoid overfitting we augment the B-spline basis with a quadratic penalty on the coefficients, typically based on the squared difference of adjacent coefficients. This can be written in the form $\lambda_j \bm{\beta}^{\top}_j \bm{S}_j \bm{\beta}_j$, where $\bm{S}_j$ is constructed based on a known difference matrix and $\lambda_j > 0$ is the smoothing parameter. A larger $\lambda_j$ forces a smoother curve, while $\lambda_j = 0$ implies no penalization. 

In this article, we will use the popular cubic B-splines and set $m_1 = \ldots = m_q = 3$. Regarding the number of knots, for simplicity we use the same number of interior knots for all $q$ covariates, i.e., $K_1 = \ldots = K_q = K$, but allow $K$ to grow with sample size. As reviewed in Section~\ref{sec:intro}, research has shown that in practice the model fit is relatively insensitive to the basis dimension provided the number of knots is large enough to offer adequate flexibility.

Let $\bm{\Psi} = (\bm{\kappa}^{\top}, \phi^{\top})^{\top}$ denote the parameters corresponding to the parametric component of the model, and $\bm{\beta} = (\bm{\beta}^{\top}_1, \ldots, \bm{\beta}^{\top}_q)^{\top}$ be the vector of $d = \sum_{j=1}^q d_j$ smoothing coefficients in the nonparametric component. The most common method of estimation and inference for GAMs is based on using the penalized log-likelihood $l_{P}(\bm{\Psi},\bm{\beta}) = \sum_{i=1}^n \ln \{f(y_i | \bm{\Psi}, \bm{\beta})\} - 2^{-1} \sum_{j=1}^q \lambda_j \bm{\beta}^{\top}_j \bm{S}_j \bm{\beta}_j$. Note the distribution of the responses $f(y_i | \bm{\Psi}, \bm{\beta})$ also depends on the covariates $\bm{x}_i$ and $\bm{u}_i$, although for ease of notation this dependence is not made explicit in the developments below. Likewise, $l_{P}(\bm{\Psi},\bm{\beta})$ depends on the smoothing parameters $\lambda_j$ but this dependence is suppressed. Given a set of $\lambda_j$'s, $l_{P}(\bm{\Psi}, \bm{\beta})$ can be maximized, for instance, via a penalized iterative reweighted least squares (P-IRLS) algorithm \citep{wood2006generalized}. Inference tools such as confidence and prediction intervals and hypothesis tests for the fitted smooths are then usually based on developing asymptotic normality theory for (a linear combination of) the $\bm{\beta}_j$'s. Finally, the penalized likelihood is also often combined with a procedure for selecting smoothing parameters based on an optimality criterion such as generalized cross-validation.

An alternative to the penalized likelihood approach is to formulate the GAM as a mixed effects model, by recognizing that the quadratic penalty in $l_{P}(\bm{\Psi},\bm{\beta})$ amounts to assuming a normal random effects distribution. For $j = 1,\ldots,q$ let $f(\bm{\beta}_j | \lambda_j) = \mathcal{N}_{d_j}\{\bm{0}, (\lambda_j \bm{S}_j)^{-1}\} = c_0\lambda_j^{2^{-1} d_j} \exp\left(-2^{-1} \lambda_j \bm{\beta}^{\top}_j \bm{S}_j \bm{\beta}_j\right)$ where $c_0$ is the normalizing constant independent of $\bm{\beta}_j$ and $\lambda_j$. Note that because of the centering constraint applied to each smooth, $\bm{S}_j$ is of full rank. We then have the marginal log-likelihood 
\begin{align} \label{eqn:marglogL}
\ell(\bm{\Psi}) = \ln \left(\int \prod\limits_{i=1}^n f(y_i | \bm{\Psi}, \bm{\beta}) \prod\limits_{j=1}^q f(\bm{\beta}_j | \lambda_j) d\bm{\beta}_j\right). 
\end{align}
Formulating the GAM as a mixed model allows us to exploit the wealth of associated mixed model methodology for inference. However, a major drawback is that for non-normal responses, equation~(\ref{eqn:marglogL}) does not possess a tractable form. One approach around this is to apply a Laplace approximation and maximize the subsequent approximated log-likelihood, $\ell_{\text{Lap}}(\bm{\Psi}) = l_P\{\bm{\Psi}, \hat{\bm{\beta}}(\bm{\Psi})\} + 2^{-1} \sum_{j=1}^q d_j\ln(\lambda_j) - 2^{-1}\ln\det\left[-\nabla^2 l_P\{\bm{\Psi}, \hat{\bm{\beta}}(\bm{\Psi})\}\right]$, where $\hat{\bm{\beta}}(\bm{\Psi})$ maximizes $l_P(\bm{\Psi}, \bm{\beta})$ for a given $\bm{\Psi}$ and the negative Hessian is with respect to $\bm{\beta}$. While the Laplace approximation approach overcomes the intractable integral, estimation and inference remain a challenge for two main reasons. First, $\bm{\beta}$ depends indirectly on $\bm{\Psi}$ as a solution of an inner optimization problem for $l_P(\bm{\Psi}, \bm{\beta})$. Second, the log determinant term on the right hand side is a non-linear function of $\bm{\Psi}$, and cannot be straightforwardly integrated into a P-IRLS algorithm. 

\section{Variational Approximations for GAMs} \label{sec:vagams}

Motivated by the challenges associated with the intractable integral in the mixed model formulation, we introduce a variational approximations framework to estimation and inference for GAMs. For the full $d$-vector of smoothing coefficients $\bm{\beta}$, consider a multivariate normal variational distribution $h(\bm{\beta} | \bm{a}, \bm{A}) = \mathcal{N}_{d}(\bm{a}, \bm{A})$ with mean vector $\bm{a} = (\bm{a}^{\top}_1, \ldots, \bm{a}^{\top}_q)^{\top}$ and a fully unstructured $d \times d$ covariance matrix $\bm{A}$ which is assumed to be positive definite. Unlike the case of independent clustered mixed models, where a block diagonal form for $\bm{A}$ was as optimal as a fully unstructured form \citep{ormerod2012gaussian}, assuming any form of block diagonality here is suboptimal. That is, if we were to constrain the variational distribution to be multivariate normal, then the optimal form for $\bm{A}$ which minimizes the Kullback-Leibler distance between the posterior and variational distributions of the smoothing coefficients is the fully unstructured one; see also the discussion below (\ref{eqn:KLdistanceminimize}). Of course, in some settings we may want to consider block diagonality for $\bm{A}$ for reasons of parsimony, since the unstructured form involves estimating $2^{-1}d(d+1)$ elements, which may be problematic when $n$ is not large compared to $d$. In Section~\ref{sec:sims}, we empirically compared unstructured and block diagonal forms for $\bm{A}$ (which is equivalent to setting $h(\bm{\beta}| \bm{a}, \bm{A}) = \prod_{j=1}^q \mathcal{N}_{d_j}(\bm{a}_j, \bm{A}_j)$ instead), and found little difference between the two empirically, but that assuming block diagonality for $\bm{A}$ actually took longer to estimate. 

Utilizing this variational distribution and applying Jensen's inequality to (\ref{eqn:marglogL}), we obtain
\begin{align*}
\ell(\bm{\Psi}) &= \ln \left(\int \frac{\prod\limits_{i=1}^n f(y_i | \bm{\Psi}, \bm{\beta}) \prod\limits_{j=1}^q f(\bm{\beta}_j | \lambda_j)}{h(\bm{\beta} | \bm{a}, \bm{A})} h(\bm{\beta} | \bm{a}, \bm{A}) d\bm{\beta}\right) \\
&\ge \int \sum\limits_{i=1}^n \ln \{f(y_i | \bm{\Psi}, \bm{\beta})\} h(\bm{\beta} | \bm{a}, \bm{A}) d\bm{\beta} +  \int \sum \limits_{j=1}^q\ln \{f(\bm{\beta}_j | \lambda_j)\} h(\bm{\beta} | \bm{a}, \bm{A}) d\bm{\beta} \nonumber \\
&\quad - \int \ln \{h(\bm{\beta} | \bm{a}, \bm{A})\} h(\bm{\beta} | \bm{a}, \bm{A}) d\bm{\beta}, \\
&= \underline{\ell}(\bm{\Psi}, \bm{\xi}) \nonumber,
\end{align*}
where $\underline{\ell}(\bm{\Psi}, \bm{\xi})$ denotes the variational log-likelihood for a GAM with additional parameters $\bm{\xi} = (\bm{a}^{\top}, \text{vech}(\bm{A})^{\top})^{\top}$, and $\text{vech}(\cdot)$ is the half-vectorization operator. Estimates of the model and variational parameters, denoted here as $(\hat{\bm{\Psi}}, \hat{\bm{\xi}})$, are then obtained by maximizing $\underline{\ell}(\bm{\Psi}, \bm{\xi})$. Note that if we let $\bm{y} = (y_1,\ldots,y_n)^{\top}$ and $\bm{\lambda} = (\lambda_1,\ldots,\lambda_q)^{\top}$, then the above inequality can be rearranged as 
\begin{align} \label{eqn:KLdistanceminimize}
0 \ge \underline{\ell}(\bm{\Psi}, \bm{\xi}) - \ell(\bm{\Psi}) = \int \ln\left\{\frac{f(\bm{\beta}|\bm{y},\bm{\Psi},\bm{\lambda})}{h(\bm{\beta} | \bm{a}, \bm{A})}\right\} h(\bm{\beta} | \bm{a}, \bm{A}) d\bm{\beta}.
\end{align}
Therefore, we see that maximizing $\underline{\ell}(\bm{\Psi}, \bm{\xi})$ is equivalent to minimizing the Kullback-Leibler distance between the posterior and variational distributions of the smoothing coefficients. The VA framework bears some resemblance to the well-known Expectation Maximization \citep[EM,][]{dempster77} algorithm for estimating models involving missing data. In many cases, the posterior distribution and hence the E-step does not present a tractable form, and so one often has to resort to Monte-Carlo methods to perform this. The VA approach can be viewed as an approximation to the EM algorithm which sacrifices exactness for computational speed: we approximate the intractable posterior distribution of the random effects by a variational distribution, allowing us to construct a fully tractable posterior expectation of the complete log-likelihood. Updating the variational parameters $\bm{\xi}$ in $\underline{\ell}(\bm{\Psi}, \bm{\xi})$ corresponds then to an approximate E-step, while updating the model parameters $\bm{\Psi}$ corresponds to an approximate M-step. 

Regardless of the distribution of the responses $y_i$, the second and third term in $\underline{\ell}(\bm{\Psi}, \bm{\xi})$ have the same form. 
\begin{align*}
\int \sum\limits_{j=1}^q \ln \{f(\bm{\beta}_j | \lambda_j)\} h(\bm{\beta} | \bm{a}, \bm{A}) d\bm{\beta} &= \frac{1}{2}\sum\limits_{j=1}^q d_j\ln(\lambda_j) - \frac{1}{2}\sum\limits_{j=1}^q \left\{\lambda_j\bm{a}_j^{\top} \bm{S}_j \bm{a}_j + \lambda_j\text{tr}(\bm{S}_j\bm{A}_j)\right\} \\
\int \ln \{h(\bm{\beta} | \bm{a}, \bm{A})\} h(\bm{\beta} | \bm{a}, \bm{A}) d\bm{\beta} &= -\frac{d}{2} \{\ln(2\pi) + 1 \}- \frac{1}{2}\ln\det(\bm{A}),
\end{align*}
where $\bm{A}_j$ is of dimension $d_j \times d_j$ and denotes the $j^{\text{th}}$ diagonal block in $\bm{A}$, and quantities related to $c_0$ have been omitted. Following some straightforward algebra then, we obtain
\begin{align} \label{eqn:valogLgeneral}
\underline{\ell}(\bm{\Psi}, \bm{\xi}) &= \int \sum\limits_{i=1}^n \ln \left\{\frac{1}{\phi}\{y_i \omega_i - b(\omega_i)\} + c(y_i,\phi)\right\} h(\bm{\beta} | \bm{a}, \bm{A}) d\bm{\beta} \nonumber \\ 
&\quad + \frac{1}{2}\sum\limits_{j=1}^q \left\{d_j \ln(\lambda_j) - \lambda_j\bm{a}_j^{\top} \bm{S}_j \bm{a}_j - \lambda_j\text{tr}(\bm{S}_j\bm{A}_j)\right\} + \frac{1}{2}\ln\det(\bm{A}).
\end{align}

One interesting feature of (\ref{eqn:valogLgeneral}) is that, conditional on all parameters, we can obtain a closed-form update for each smoothing parameter, $\hat{\lambda}_j = d_j(\text{tr}(\bm{S}_j\bm{A}_j) + \bm{a}^{\top}_j\bm{S}_j\bm{a}_j)^{-1}$. This is particularly appealing as it implies the VA approach simultaneously provides estimates of the smoothing coefficients and controls the degree of penalization on these coefficients. We elaborate on this critical aspect of inference in Section~\ref{subsec:smoothingparam}.

We now focus on three particular responses of interest where GAMs are commonly applied.

\subsection{Poisson Responses} \label{subsec:pois}

For Poisson GAMs, we have $\phi = 1$, $\omega_i = \ln(\mu_i) = \bm{x}_i^{\top}\bm{\kappa} + \sum_{j=1}^q \bm{z}^{\top}_{ij}\bm{\beta}_j$, and $b(\omega_i) = \exp(\mu_i)$. Let $\bm{z}_i = (\bm{z}^{\top}_{i1},\ldots,\bm{z}^{\top}_{iq})^{\top}$ such that $\sum_{j=1}^q \bm{z}^{\top}_{ij}\bm{a}_j = \bm{z}^{\top}_i \bm{a}$. Then it is straightforward to show 
\begin{align*} 
\underline{\ell}_{\text{Pois}}(\bm{\Psi}, \bm{\xi}) &= \sum\limits_{i=1}^n \left\{y_i\left(\bm{x}_i^{\top}\bm{\kappa} + \bm{z}^{\top}_{i}\bm{a}\right) - \exp\left(\bm{x}_i^{\top}\bm{\kappa} + \bm{z}^{\top}_{i}\bm{a} + \frac{1}{2} \bm{z}^{\top}_{i} \bm{A} \bm{z}_{i}\right)\right\}\\
&\quad + \frac{1}{2}\sum\limits_{j=1}^q \left\{d_j \ln(\lambda_j) - \lambda_j\bm{a}_j^{\top} \bm{S}_j \bm{a}_j - \lambda_j\text{tr}(\bm{S}_j\bm{A}_j)\right\} + \frac{1}{2}\ln\det(\bm{A}).
\end{align*}
The variational likelihood is fully closed form and can be easily maximized by iterative updating the model and variational parameters (for a given set of $\lambda_j$'s). For instance, we can iterate between the following steps. 
\begin{enumerate}
 \item Update all coefficients $(\bm{\kappa}, \bm{a})$ by fitting a log-link Poisson GLM with linear predictor $\bm{x}^\top_i\bm{\kappa} + \bm{z}^{\top}_{i}\bm{a}$, an offset equal to $2^{-1} \bm{z}^{\top}_{i} \bm{A} \bm{z}_{i}$, and a quadratic penalty of $2^{-1}\sum_{j=1}^q \lambda_j\bm{a}_j^{\top} \bm{S}_j \bm{a}_j$. The score equations for these are provided in the Supplementary Material.
 \item The score equation for the covariance matrix is
 \begin{align*}
 \frac{\partial \underline{\ell}_{\text{Pois}}(\bm{\Psi}, \bm{\xi})}{\partial \bm{A}} &= -\frac{1}{2}\sum\limits_{i=1}^n \bm{z}_{i}\bm{z}^{\top}_{i} \exp\left(\bm{x}_i^{\top}\bm{\kappa} + \bm{z}^{\top}_{i}\bm{a} + \frac{1}{2} \bm{z}^{\top}_{i} \bm{A} \bm{z}_{i}\right) + \frac{1}{2}\left(\bm{A}^{-1} - \bm{S}_{\bm{\lambda}}\right),
 \end{align*}
 where $\bm{S}_{\bm{\lambda}}$ is a $d \times d$ block diagonal matrix formed by taking blocks $\lambda_j\bm{S}_j$ for $j = 1,\ldots,q$. The above could then be used directly as part of a Quasi-Newton optimization routine, subject to $\bm{A}$ being positive definite. In settings when $d$ is large, it may be numerically more stable to parameterize and hence optimize with respect to the Cholesky decomposition of $\bm{A}$. Alternatively, solving for the score equation suggests a fixed point iterative estimator $\bm{A}^{(1)} = \left\{\bm{S}_{\bm{\lambda}} + \sum_{i=1}^n \bm{z}_{i}\bm{z}^{\top}_{i} \exp\left(\bm{x}_i^{\top}\bm{\kappa} + \bm{z}^{\top}_{i}\bm{a} + 2^{-1} \bm{z}^{\top}_{i} \bm{A}^{(0)} \bm{z}_{i}\right)\right\}^{-1}$ where $\bm{A}^{(1)}$ and $\bm{A}^{(0)}$ are the new and current estimate, respectively. 
\end{enumerate}

\subsection{Normal Responses} \label{subsec:normal}

Assuming an identity link, let $\tilde{r}_i = y_i - \bm{x}_i^{\top}\bm{\kappa}$ denote the partial residual from the parametric component. Then $\ln \{f(y_i | \bm{\Psi}, \bm{\beta})\} = -2^{-1}\ln(2\pi\phi) - (2\phi)^{-1} \left\{\tilde{r}_i^2 + \bm{\beta}^{\top} \bm{z}_{i}\bm{z}^{\top}_{i} \bm{\beta}\right\} + \phi^{-1} \tilde{r}_i \bm{z}^{\top}_{i} \bm{\beta}$ and we obtain
\begin{align*} 
\underline{\ell}_{\text{Norm}}(\bm{\Psi}, \bm{\xi}) &= -\frac{n}{2}\ln(\phi) -\frac{1}{2\phi} \sum\limits_{i=1}^n \left(y_i - \bm{x}_i^{\top}\bm{\kappa} - \bm{z}^{\top}_{i}\bm{a} \right)^2 - \frac{1}{2\phi}\sum\limits_{i=1}^n \bm{z}^{\top}_{i} \bm{A}\bm{z}_{i} \\
&\quad + \frac{1}{2}\sum\limits_{j=1}^q \left\{d_j \ln(\lambda_j) - \lambda_j\bm{a}_j^{\top} \bm{S}_j \bm{a}_j - \lambda_j\text{tr}(\bm{S}_j\bm{A}_j)\right\} + \frac{1}{2}\ln\det(\bm{A}),
\end{align*}
where constant terms with respect to $(\bm{\Psi}, \bm{\xi})$ have been omitted. Once more, the variational likelihood has a fully closed-form and can be easily maximized as follows.
\begin{enumerate}
 \item Update all coefficients $(\bm{\kappa}, \bm{a})$ by fitting a linear model with linear predictor $\bm{x}^{\top}_i\bm{\kappa} + \bm{z}^{\top}_{i}\bm{a}$, and a quadratic penalty of $2^{-1}\sum_{j=1}^q \lambda_j\bm{a}_j^{\top} \bm{S}_j \bm{a}_j$. The score equations for these are provided in the Supplementary Material.
 \item Update the dispersion parameter as $\phi = n^{-1} \sum_{i=1}^n \left\{(\tilde{r}_i - \bm{z}^{\top}_{i}\bm{a})^2 + \bm{z}_{i}^{\top} \bm{A} \bm{z}_{i}\right\}$. 
 \item Solve the score equation for the covariance matrices to obtain the closed-form update $\bm{A} = \left(\bm{S}_{\bm{\lambda}} + \phi^{-1}\sum_{i=1}^n \bm{z}_{i}\bm{z}^{\top}_{i}\right)^{-1}$.
\end{enumerate}
Note step 1 is independent of $\bm{A}$ and $\phi$. In the case of normal responses, the VA approach coincides exactly with using the EM algorithm to fit the GAM. In particular, for GAMs it is straightforward to show the posterior distribution of $\bm{\beta}$ is multivariate normal with covariance $\left(\bm{S}_{\bm{\lambda}} + \phi^{-1}\sum_{i=1}^n \bm{z}_{i}\bm{z}^{\top}_{i}\right)^{-1}$ and mean vector $\left(\bm{S}_{\bm{\lambda}} + \phi^{-1}\sum_{i=1}^n \bm{z}_{i}\bm{z}^{\top}_{i}\right)^{-1}\left(\sum_{i=1}^n (y_i - \bm{x}^{\top}_i\bm{\kappa})\bm{z}_i \right)$, which corresponds exactly to the formulas for $\bm{A}$ and $\bm{a}$, respectively.

\subsection{Bernoulli Responses} \label{subsec:bern}

A number of approaches have been proposed for handling Bernoulli responses in the VA framework. \citet{ormerod2012gaussian} considered a logit link and accepted the fact that this resulted in a term, $b(\omega) = \ln\{1 + \exp(\eta)\}$, whose expectation with respect to the normal variational distribution $h(\bm{\beta} | \bm{a}, \bm{A})$ would require (univariate) Monte Carlo integration. In an attempt to get around this problem, \citet{hui2016variational} considered the probit link and exploited the fact that the model could be reparametrized by introducing a normally distributed auxiliary variable. However, a downside with this approach, which was not picked up by \citet{hui2016variational}, is that the estimate of $\bm{A}$ can be quite biased and variability of the posterior distribution of $\bm{\beta}$ (as approximated by the variational distribution) tends to be underestimated. Following on from \citet{blei07}, we propose an alternative method of handling Bernoulli responses for the VA that is both simple and overcomes the problems mentioned above with other approaches. By using the canonical logit link, it is straightforward to show that the variational likelihood for a Bernoulli GAM is given by $\underline{\ell}(\bm{\Psi}, \bm{\xi}) = \sum_{i=1}^n \left[y_i\left(\bm{x}^{\top}_i\bm{\kappa} + \bm{z}^{\top}_{i} \bm{a}\right) - \int \ln\left\{1+\exp\left(\bm{x}^{\top}_i\bm{\kappa} + \bm{z}^{\top}_{i} \bm{\beta}\right)\right\} h(\bm{\beta} | \bm{a}, \bm{A}) d\bm{\beta} \right] \\ + 2^{-1}\sum_{j=1}^q \left\{d_j \ln(\lambda_j) - \lambda_j\bm{a}_j^{\top} \bm{S}_j \bm{a}_j - \lambda_j\text{tr}(\bm{S}_j\bm{A}_j)\right\} + 2^{-1}\ln\det(\bm{A})$, where constant terms with respect to $(\bm{\Psi}, \bm{\xi})$ have been omitted. We immediately see the problem of the intractable integral. However, by exploiting Jensen's inequality once more, we obtain 
\begin{align*}  
\underline{\ell}(\bm{\Psi}, \bm{\xi}) &\ge \sum\limits_{i=1}^n \left[y_i\left(\bm{x}^{\top}_i\bm{\kappa} + \bm{z}^{\top}_{i} \bm{a}\right) - \ln\left\{1+\exp\left(\bm{x}^{\top}_i\bm{\kappa} + \bm{z}^{\top}_{i} \bm{a} + \frac{1}{2}\bm{z}^{\top}_i\bm{A}\bm{z}_i\right)\right\} \right] \\
&\quad + \frac{1}{2}\sum\limits_{j=1}^q \left\{d_j \ln(\lambda_j) - \lambda_j\bm{a}_j^{\top} \bm{S}_j \bm{a}_j - \lambda_j\text{tr}(\bm{S}_j\bm{A}_j)\right\} + \frac{1}{2}\ln\det(\bm{A}) \\
&\triangleq \underline{\ell}_{\text{Bern}}(\bm{\Psi}, \bm{\xi}).
\end{align*}	
The above objective function is now fully tractable and can be maximized in a similar way to the Poisson response case with two differences: in step 1, we fit a logistic regression model with linear predictor $\bm{x}^{\top}_i\bm{\kappa} + \bm{z}^{\top}_{i}\bm{a}$, an offset equal to $2^{-1} \bm{z}^{\top}_{i} \bm{A} \bm{z}_{i}$, and a quadratic penalty of $2^{-1}\sum_{j=1}^q \lambda_j\bm{a}_j^{\top} \bm{S}_j \bm{a}_j$, and in step 2 the fixed point iterative estimator becomes $\bm{A}^{(1)} = \left(\bm{S}_{\bm{\lambda}} + \sum_{i=1}^n w^{(0)}_i \bm{z}_{i}\bm{z}^{\top}_{i} \right)^{-1}$ where $w^{(0)}_i = \left\{1 + \exp\left(-\bm{x}_i^{\top}\bm{\kappa} - \bm{z}^{\top}_{i}\bm{a} - 2^{-1} \bm{z}^{\top}_{i} \bm{A}^{(0)} \bm{z}_{i}\right)\right\}^{-1}$. 

While it remains a lower bound to the marginal log-likelihood, the second application of Jensen's inequality means $\underline{\ell}_{\text{Bern}}(\bm{\Psi}, \bm{\xi})$ is now suboptimal in the sense that we are no longer minimizing the Kullback-Leibler distance between the posterior and variational distributions of the smoothing coefficients. However, it turns out that this additional approximation does not have a detrimental effect on the asymptotic performance of the estimates, and the VA framework here still produces consistent and asymptotically normal estimates; see Section~\ref{sec:sims} and also \citet{knowles11} for a discussion of related bounds in binary and multinomial regression. 


\section{Inference} \label{sec:inference}

Equation~\eqref{eqn:KLdistanceminimize} showed that maximizing the variational log-likelihood is equivalent to minimizing the Kullback-Leibler distance between the true posterior distribution $f(\bm{\beta}|\bm{y},\bm{\Psi},\bm{\lambda})$ and the variational distribution $h(\bm{\beta} | \bm{a}, \bm{A})$ of the smoothing coefficients. Along with the normality assumption on $h(\bm{\beta} | \bm{a}, \bm{A})$, this suggests that inference on the smooth component of the GAM can be obtained directly from the estimation process. For example, $\hat{\bm{a}}$ serves as the variational version of both the empirical Bayes and maximum a-posteriori estimate of the smoothing coefficients, while $\hat{\bm{A}}$ is an estimate of the posterior covariance matrix. The multivariate normality of $h(\bm{\beta} | \bm{a}, \bm{A})$ also means that we can easily construct confidence intervals and hypothesis tests. For example, we can test the null hypothesis $H_0: \bm{\beta}_j = \bm{0}$ by comparing the Wald statistic $W = \hat{\bm{a}}_j^{\top}(\hat{\bm{A}}^{-1})_{j}\bm{a}_j$ to a Chi-squared distribution with $d_j$ degrees of freedom, where $(\hat{\bm{A}}^{-1})_{j}$ denotes the diagonal submatrix block of $\hat{\bm{A}}^{-1}$ relating to covariate $j = 1,\ldots,q$. This idea extends naturally to inference for the estimated smooth functions. On the linear predictor scale, the fitted smooth for covariate $j = 1,\ldots,q$ at observation $i = 1,\ldots,n$ is given by $\bm{z}^{\top}_{ij}\hat{\bm{a}}_j$, and a $(1-\alpha)100\%$ pointwise confidence interval is then given by $\{\bm{z}^{\top}_{ij}\hat{\bm{a}}_j \pm \Phi^{-1}_{1-\alpha/2}(\bm{z}^{\top}_{ij}\hat{\bm{A}}_j\bm{z}_{ij})^{1/2}\}$ where $\Phi^{-1}_{1-\alpha/2}$ is the $(1-\alpha/2)$-th quantile of the standard normal distribution. For simultaneous confidence bands, one straightforward way to construct this within the VA framework would be as follows: for covariate $j = 1,\ldots,q$, simulate a large number of realizations $g = 1,\ldots,G$ from the variational distribution $h(\bm{\beta} | \hat{\bm{a}}, \hat{\bm{A}}) = \mathcal{N}_d(\hat{\bm{a}}, \hat{\bm{A}})$. Then, for a grid of $l = 1,\ldots,M$ values spanning the range of covariate $j$, calculate the corresponding values of $C_g = \max_{l = 1,\ldots,M}  \left\{\bm{z}^{\top}_{(l)j}\hat{\bm{a}}_j\left(\bm{z}^{\top}_{(l)j}\hat{\bm{A}}_j\bm{z}_{(l)j}\right)^{-1/2}\right\}$, where $\bm{z}_{(l)j}$ is the vector of basis function values for covariate $j$, evaluated at the $l$-th grid value. By empirically determining the $(1-\alpha/2)$ quantile of the values $\{C_g; g = 1,\ldots,G\}$, which we denote as $c_{1-\alpha/2}$, a $(1-\alpha)100\%$ simultaneous confidence band at observation $i = 1,\ldots,n$ is then given by $\{\bm{z}^{\top}_{ij}\hat{\bm{a}}_j \pm c_{1-\alpha/2} (\bm{z}^{\top}_{ij}\hat{\bm{A}}_j\bm{z}_{ij})^{1/2}\}$. We do acknowledge however that construction of simultaneous confidence bands remains an active area of research in semiparametric regression \citep[e.g.,][]{krivobokova2010}, as different bands can be obtained depending on the formulation of the model to adopt. 

\subsection{Parametric Component} \label{subsec:parametricbit}

Previously, to obtain standard errors for the fixed effect parameters of the model, \citet{ormerod2012gaussian} and \citet{hui2016variational} proposed obtaining an information matrix by directly calculating the negative Hessian of the variational log-likelihood. While this seems logical at first glance, one downside of this approach is that it involves deriving a Hessian with respect to both model \emph{and} variational parameters, i.e., both $\bm{\Psi}$ and $\bm{\xi}$. In general, this can be computationally intensive given the increased number of parameters appearing in the Hessian, although with independent cluster models \citep[as is the case in][]{ormerod2012gaussian,hui2016variational} there are certain structures in the Hessian matrix that can be exploited. Another, more conceptual drawback is that it seems superfluous to take derivatives with respect to variational parameters, given our objective here is to quantify the uncertainty only in the parameters that define the GAM, i.e., $\bm{\Psi}$ and possibly the smoothing parameters $\bm{\lambda}$.

We propose an alternative method of obtaining standard errors for the parametric component, based on exploiting Louis' method for obtaining the observed information matrix \citep{louis82}. For the mixed model formulation of a GAM given in (\ref{eqn:marglogL}), write the complete log-likelihood as $\ell_{\text{com}}(\bm{\Psi},\bm{\beta}) = \sum_{i=1}^n \ln\{f(y_i | \bm{\Psi}, \bm{\beta})\} + 2^{-1}\sum_{j=1}^q d_j \ln(\lambda_j) - 2^{-1} \sum_{j=1}^q \lambda_j \bm{\beta}^{\top}_j \bm{S}_j\bm{\beta}_j$, ignoring constants. By Louis' method, the observation information matrix of the marginal log-likelihood, $\bm{I}(\bm{\Psi},\bm{\lambda}) = -\partial^2 \ell(\bm{\Psi})/\partial (\bm{\Psi},\bm{\lambda}) \partial (\bm{\Psi},\bm{\lambda})^{\top}$ can be calculated as
\begin{align*}
\bm{I}(\bm{\Psi},\bm{\lambda}) =& \ \int -\frac{\partial^2 \ell_{\text{com}}(\bm{\Psi},\bm{\beta})}{\partial (\bm{\Psi},\bm{\lambda}) \partial (\bm{\Psi},\bm{\lambda})^{\top}} f(\bm{\beta}|\bm{y},\bm{\Psi},\bm{\lambda}) \: d\bm{\beta} \\
&- \int \left(\frac{\partial \ell_{\text{com}}(\bm{\Psi},\bm{\beta})}{\partial (\bm{\Psi},\bm{\lambda})}\right)\left(\frac{\partial \ell_{\text{com}}(\bm{\Psi},\bm{\beta})}{\partial (\bm{\Psi},\bm{\lambda})}\right)^{\top} f(\bm{\beta}|\bm{y},\bm{\Psi},\bm{\lambda}) \: d\bm{\beta} + \left(\frac{\partial \ell(\bm{\Psi})}{\partial (\bm{\Psi},\bm{\lambda})}\right) \left(\frac{\partial \ell(\bm{\Psi})}{\partial (\bm{\Psi},\bm{\lambda})}\right)^{\top}.
\end{align*}
Note the first and second derivatives of $\ell_{\text{com}}(\bm{\Psi},\bm{\beta})$ with respect to $(\bm{\Psi},\bm{\lambda})$ are comparably straightforward to compute. Within the VA framework, we propose modifying the above formula in two ways: I) Replace $f(\bm{\beta}|\bm{y}, \bm{\Psi},\bm{\lambda})$ with the variational distribution $h(\bm{\beta} | \bm{a}, \bm{A}) = \mathcal{N}_d(\bm{a}, \bm{A})$. As in (\ref{eqn:KLdistanceminimize}), maximizing the VA likelihood is equivalent to minimizing the Kullback-Leibler distance between the true posterior and variational distribution, and so this is a logical substitution to make; II) Assume that at the VA estimates, the score equation for the marginal log-likelihood is approximately zero, $\partial \ell(\hat{\bm{\Psi}}) /\partial (\bm{\Psi},\bm{\lambda}) \approx \bm{0}$. In standard maximum likelihood estimation, this is exactly equal to the zero vector by definition. While this is not guaranteed to be true also for the VA estimates, except for the normal response case (see Section~\ref{subsec:normal}), it should nevertheless be close to zero and negligible compared to the first two terms in Louis' method given the estimates are maximizing a lower bound to the marginal log-likelihood. 
Applying the two modifications above leads us to a new, variational observed information matrix evaluated at $(\hat{\bm{\Psi}},\hat{\bm{\lambda}})$
\begin{align} \label{eqn:infomat}
\bm{I}_{\text{v}}(\hat{\bm{\Psi}}, \hat{\bm{\lambda}}) &= \left[\int -\frac{\partial^2 \ell_{\text{com}}(\bm{\Psi},\bm{\beta})}{\partial (\bm{\Psi},\bm{\lambda}) \partial (\bm{\Psi},\bm{\lambda})^{\top}} h(\bm{\beta}|\hat{\bm{a}},\hat{\bm{A}}) \: d\bm{\beta}\right]_{(\hat{\bm{\Psi}},\hat{\bm{\lambda}})} \nonumber \\
&\quad - \left[\int \left(\frac{\partial \ell_{\text{com}}(\bm{\Psi},\bm{\beta})}{\partial (\bm{\Psi},\bm{\lambda})}\right)\left(\frac{\partial \ell_{\text{com}}(\bm{\Psi},\bm{\beta})}{\partial (\bm{\Psi},\bm{\lambda})}\right)^{\top} h(\bm{\beta}|\hat{\bm{a}},\hat{\bm{A}}) \: d\bm{\beta}\right]_{(\hat{\bm{\Psi}},\hat{\bm{\lambda}})}.
\end{align}
As usual, the above can be inverted to obtain standard errors and hence be used as the basis for confidence intervals and Wald tests for the parametric coefficients $\bm{\kappa}$. 
Note that $h(\bm{\beta} | \bm{a}, \bm{A}) = \mathcal{N}_d(\bm{a}, \bm{A})$ is easy to sample from and the derivatives of the complete log-likelihood $\ell_{\text{com}}(\bm{\Psi},\bm{\beta})$ are easy to calculate. Thus the variational information matrix $\bm{I}_{\text{v}}(\hat{\bm{\Psi}}, \hat{\bm{\lambda}})$ is straightforward and comparably efficient to construct. In the Supplementary Material, we provide more information regarding the calculation of (\ref{eqn:infomat}) for the three responses of interest in this article. 
Note also that if $\bm{\lambda}$ is estimated by maximizing the variational likelihood, then (\ref{eqn:infomat}) could in principle be used to obtain standard errors for the smoothing parameters. 

We conclude this section by emphasizing that the above covers only some of the primary inferences one may wish to perform in semiparametric regression. Other aspects, such as variable selection, the choice of the number of knots, 
and how they can be implemented under the VA framework, are avenues of future research. One of the most critical aspects of inference, however, that does warrant further attention is the choice of smoothing parameters, and we address this in the next section.

\subsection{Choosing the Smoothing Parameters} \label{subsec:smoothingparam}

As mentioned at the end of Section~\ref{sec:vagams}, one of the appealing aspects of the VA framework for GAMs is that it provides closed-form updates of the smoothing parameters, $\hat{\lambda}_j = d_j(\text{tr}(\bm{S}_j\bm{A}_j) + \bm{a}^{\top}_j\bm{S}_j\bm{a}_j)^{-1}$ for $j = 1,\ldots,q$. This formula can be naturally integrated into the sequence of update steps outlined for the three responses discussed in Sections~\ref{subsec:pois}-\ref{subsec:bern}. This approach to selecting the smoothing parameters is quite advantageous: not only is its closed form easily calculated, but it also means we do not have to employ a separate method for choosing $\bm{\lambda}$ that is external to the estimation of the model parameters.
It turns out that this simple formula for choosing the smoothing parameters in the VA framework bears a close resemblance to estimating smoothing parameters via the EM algorithm and the Laplace approximation. For the former, recall the discussion in Section~\ref{sec:vagams} that updating the variational parameters and model parameters in $\underline{\ell}(\bm{\Psi}, \bm{\xi})$ corresponds roughly to an approximate E-step and M-step respectively. On the other hand, under the mixed effects parameterization for a GAM, the smoothing parameters correspond simply to a set of inverse variance (precision) components, and thus we could treat $\bm{\lambda}$ as model parameters and update them as part of the M-step of the EM algorithm, in addition to $\bm{\Psi}$. The closed-form update $\hat{\lambda}_j$ above is nothing more than the analog of this in the VA framework, i.e., an approximate M-step for updating the smoothing parameters. Turning to the Laplace approximation, if the Laplace approximated log-likelihood $\ell_{\text{Lap}}(\bm{\Psi})$ in Section~\ref{sec:gams} is treated as a function of both $\bm{\Psi}$ \emph{and} $\bm{\lambda}$, then solving $\partial \ell_{\text{Lap}}(\bm{\Psi})/\partial \lambda_j = 0$ leads to the formula $\hat{\lambda}_{j,\text{Lap}} = d_j\left(\text{tr}\left(\bm{S}_j\bm{G}\{\bm{\Psi}, \hat{\bm{\beta}}(\bm{\Psi})\}_j\right) + \hat{\bm{\beta}}(\bm{\Psi})^{\top}_j \bm{S}_j \hat{\bm{\beta}}(\bm{\Psi})_j\right)^{-1}$ where $\bm{G}\{\bm{\Psi}, \hat{\bm{\beta}}(\bm{\Psi})\}$ is the $d_j \times d_j$ submatrix of $[-\nabla^2 l_P\{\bm{\Psi}, \hat{\bm{\beta}}(\bm{\Psi})\}]^{-1}$ associated with covariate $j = 1,\ldots,q$. Note that the forms of $\hat{\lambda}_{j,\text{Lap}}$ and $\hat{\lambda}_{j}$ are similar: the Laplace approximation estimate of the smoothing parameter depends on $\hat{\bm{\beta}}(\bm{\Psi})$ and $[-\nabla^2 l_P\{\bm{\Psi}, \hat{\bm{\beta}}(\bm{\Psi})\}]^{-1}$ which, assuming the joint likelihood function $\prod_{i=1}^n f(y_i | \bm{\Psi}, \bm{\beta}) \prod_{j=1}^q f(\bm{\beta}_j | \lambda_j)$ for fixed $\bm{\Psi}$ is approximately normally distributed, can be interpreted as the mean and covariance matrix of this normal distribution. In contrast, the VA estimate of the smoothing parameter depends on $\bm{a}$ and $\bm{A}$, which are the mean and covariance matrix for the normal distribution that best approximates the posterior distribution of $\bm{\beta}$ in the Kullback-Leibler sense. Therefore, we see that VA and the Laplace approximation produce almost identical formulas for estimating the smoothing parameters, with the differences between the two methods arising as a natural consequence of where the normality assumption is made, i.e., the joint likelihood function versus the posterior distribution; see the discussion following equation~(\ref{eqn:KLdistanceminimize}).



\section{Asymptotic Theory} \label{sec:asymptotic}

We study the asymptotic properties of the VA approach to GAMs under a similar setting to that of \citet{yu2002penalized} and \citet{lu2008penalized} among others. That is, for a set of $n$ independent observations we assume the true regression model takes the form $g(\mu_i) = \eta_i = \bm{x}^{\top}_i \bm{\kappa}^0 + \sum_{j=1}^q \bm{z}^{\top}_{ij}\bm{\beta}^0_j$ for $i  = 1,\ldots,n$, where $\bm{\kappa}^0$ and $\bm{\beta}^0 = \{(\bm{\beta}^{0}_1)^{\top}, \ldots, (\bm{\beta}^{0}_q)^{\top}\}^{\top}$ are the true regression coefficients for the parametric and smooth components respectively, and develop consistency and asymptotic normality of the VA estimates. Unlike \citet{yu2002penalized} and \citet{lu2008penalized} however, where the number of basis functions was fixed, we allow $d = \dim(\bm{\beta})$ to grow polynomially with sample size. As reviewed in Section~\ref{sec:intro}, existing literature has shown that growing the number of basis functions polynomially is sufficient to accurately capture the true mean response. 

Let $\bm{\Psi}^0 = \{(\bm{\kappa}^{0})^{\top}, (\phi^{0})^{\top}\}^{\top}$ denote the parameters in the parametric component of the GAM (including the dispersion parameter) and $\bm{\theta} = (\bm{\Psi}^{\top},\bm{\beta}^{\top})^{\top}$ the full parameter vector in the GAM. Then we let $\theta^0 = \{(\bm{\Psi}^{0})^{\top},(\bm{\beta}^{0})^{\top}\}^{\top}$ denote the true parameter point, $\hat{\bm{\theta}} = (\hat{\bm{\Psi}}^{\top},\hat{\bm{a}}^{\top})^{\top}$ denote the VA estimates for a given smoothing parameter. We require the following regularity assumptions:
\begin{itemize}
 \item[(C1)] For all $i = 1,\ldots,n$, the probability density $f(y_i | \bm{\Psi}, \bm{\beta}) = f(y_i | \bm{\theta})$ has common support and is at least three times differentiable in $\bm{\theta}$. Furthermore, the model is identifiable in $\bm{\theta}$, i.e., if $\bm{\theta}' \ne \bm{\theta}$ then $f(y_i | \bm{\theta}') \ne f(y_i | \bm{\theta})$. 
 \item[(C2)] For all $i = 1,\ldots,n$ there exists a constant $C_1$ such that $\|(\bm{x}^{\top}_{i},\bm{z}^{\top}_{i})^{\top}\|_{\infty} < C_1 < \infty$ where $\|\cdot\|_{\infty}$ is the infinity norm.
 \item[(C3)] The true parameter point $\bm{\theta}^0$ is in the interior of the parameter space $\bm{\Theta}$, and satisfies $\text{s}(\bm{\theta}) = \text{E}\left(\partial \ln \{f(y_1 | \bm{\theta})\}/\partial\bm{\theta}\right)$ and $\text{s}(\bm{\theta}^0) = \bm{0}$. Furthermore, for all $i = 1,\ldots,n$ there exists a constant $C_2$ such that the linear predictor satisfies $|\eta_i| < C_2 < \infty$ at $\bm{\theta}^0$.
 \item[(C4)] The Fisher information matrix $\mathcal{I}(\bm{\theta}) = \text{E}\left(-\partial^2 \ln \{f(y_1 | \bm{\theta})\}/\partial\bm{\theta}\partial \bm{\theta}^{\top}\right)$ is finite and positive definite at $\bm{\theta}^0$, with a minimum eigenvalue that is bounded away from zero. 
 \item[(C5)] There exists an open subset $\bm{\Theta}_s \in \bm{\Theta}$ containing $\bm{\theta}^0$ such that for all $r,s,t = 1,\ldots,\dim(\bm{\theta})$, there exist functions $F_{rs}(y_1|\bm{\theta})$ and $G_{rst}(y_1|\bm{\theta})$ satisfying $|\partial^2 \ln \{f(y_1 | \bm{\theta})\}/\partial\theta_r \partial\theta_s| \le F_{rs}(y_1|\bm{\theta})$ and $|\partial^3 \ln \{f(y_1 | \bm{\theta})\}/\partial\theta_r \partial\theta_s \partial\theta_t| \le G_{rst}(y_1|\bm{\theta})$ for all $\theta \in \bm{\Theta}_s$. Furthermore, there exists constants $C_3$ and $C_4$ such that $\text{E}\{F_{rs}^2(y_1|\bm{\theta})\} < C_3 < \infty$  and $\text{E}\{G_{rst}^2(y_1|\bm{\theta})\} < C_4 < \infty$ for all $r,s,t = 1,\ldots,\dim(\bm{\theta})$.
 \item[(C6)] $d = o(n^{1/4})$; \hspace{2cm} (C6') $d = o(n^{1/5})$.
\end{itemize}
Conditions (C1)-(C5) resemble regularity conditions often employed for studying maximum likelihood estimation in a variety of regression models \citep[e.g.,][]{fan04,hui16crepe}. This should not be surprising since, as we shall see below, if the effect of the smoothing parameter is asymptotically negligible then the GAM resembles a GLM with the relevant B-spline basis functions, i.e., with covariate vector $(\bm{x}^{\top}_{i},\bm{z}^{\top}_{i})^{\top}$. Condition (C2) implies
that the covariates are non-stochastic, which is done primarily to simplify the proofs, and could be relaxed to permit $(\bm{x}^{\top}_{i},\bm{z}^{\top}_{i})^\top$ to be random but stochastically bounded although we do not pursue such an extension here. Conditions (C6) and (C6') specify the polynomial rate of growth for the spline bases dimension in the GAM, and both are milder than what is recommended in the literature, e.g., for cubic P-splines, \citet{kauermann2009some} required $d = O(n^{1/9})$. 

We first establish the following result concerning the behavior of the variational estimator $\hat{\bm{A}}$.
\begin{lemma} \label{lem:Abehave}
Under Conditions (C1)-(C4), and if $\lambda_j = o(n^{1/2})$ for all $j = 1,\ldots,q$, then for any $\bm{\theta}$ satisfying $\|\bm{\theta} - \bm{\theta}^0\| = O(d^{1/2}n^{-1/2})$ it holds that $\hat{\bm{A}} = O_p(n^{-1})$ element-wise. 
\end{lemma}
The above result is not overly surprising since $\bm{A}$ resembles a posterior covariance matrix for the smoothing coefficients, and like an inverse information matrix we would expect this to grow at the rate of the sample size. 

\begin{theorem} \label{thm:estconsist}
Under Conditions (C1)-(C6), and if $\lambda_j = o(n^{1/2})$ for all $j = 1,\ldots,q$, then the VA estimates satisfy $\left\|\hat{\bm{\theta}} - \bm{\theta}^0\right\| = O_p(d^{1/2}n^{-1/2})$. 
\end{theorem}
The above result establishes estimation consistency of the VA estimates i.e., if the effect of the smoothing parameters is asymptotically negligible, then the estimates converge to the true parameter values. This implies that we can consistently estimate both the parametric and semiparametric (non-linear) components in the model. This notion that the amount of the smoothing becomes asymptotically negligible is encapsulated in the assumption $\lambda_j = o(n^{1/2})$, which is akin to the assumption made regarding the smoothing parameters in \citet{yu2002penalized} and \citet{lu2008penalized}. Specifically, it implies that the random effects distribution for the smoothing coefficients is asymptotically dominated by the likelihood of the responses, because the inverse variance (precision) components are growing at a smaller rate than $n$. Since we are studying asymptotics under the framework that the true model is a GLM with covariates $(\bm{x}^{\top}_{i},\bm{z}^{\top}_{i})^{\top}$, then as the random effects distribution (or penalization due to smoothing) becomes asymptotically negligible, the GAM ``tends to'' the true GLM and consistency can be obtained from this. 

Note also the rate of convergence in Theorem~\ref{thm:estconsist} is the same as many (penalized) maximum likelihood estimators in regression models with a growing number of covariates \citep{fan04}, even though the estimates are based on an \emph{approximation} to the likelihood. This is an interesting result that stems from the fact that the smoothing coefficients, and in turn, the variational distribution, in a GAM is estimated from all $n$ observations.Notably, it contrasts to other applications of VA such as generalized linear mixed models and latent variable models, where the variational distribution is estimated only from a fraction of the total sample size, and therefore a price is paid in terms of convergence rate \citep{ormerod2012gaussian,hui2016variational}. We also remark that if the number of basis functions $d$ was fixed, then the familiar $n^{1/2}$-consistency arises; the slower convergence rate here is brought about due to the polynomial rate of growth of $d$ under Condition (C6)-(C6'). 

\begin{theorem} \label{thm:asympnorm}
Under Conditions (C1)-(C5) and (C6'), and if $\lambda_j = o(d^{-1}n^{1/2})$ for all $j = 1,\ldots,q$, then the VA estimates of the parametric component in a GAM also satisfies $n^{1/2}\left(\hat{\bm{\kappa}} - \bm{\kappa}^0\right) \xrightarrow{d} \mathcal{N}\left(\bm{0}, \mathcal{I}^{-1}(\bm{\theta}^0)_{\bm{\kappa}}\right)$, where $\mathcal{I}^{-1}(\bm{\theta}^0)_{\bm{\kappa}}$ is $p \times p$ submatrix of $\mathcal{I}^{-1}(\bm{\theta}^0)$ associated with $\bm{\kappa}$. 
\end{theorem}
The above result establishes the corresponding asymptotic normality of the VA estimates, i.e., if the effect of the smoothing parameters is asymptotically negligible as outlined above, then the estimates of the parametric component are asymptotically normally distributed with the covariance matrix equal to the relevant portion of the inverse Fisher information matrix in Condition (C4). 

\section{Simulation Study} \label{sec:sims} 

We compared the VA framework to several estimation and inference methods currently available for GAMs. Specifically, we considered the following methods: 1) VA with an unstructured form for $\bm{A}$ (VA-Unstruc); 2) A penalized likelihood approach using \texttt{mgcv}, with all settings set at the default options (mgcv-Default); 3) \texttt{mgcv} using P-splines and all other settings set at the default (mgcv-P--splines); 4) A mixed model approach using \texttt{gamm4} \citep{gamm4} with P-splines and all other settings set at the default (gamm4). Note that methods 1 and 4 employ a mixed model framework for GAMs, while methods 2 and 3 employ a penalized likelihood framework. It is worth pointing out that the \texttt{gamm4} package is designed more for fitting generalized additive mixed models, although it is used here to fit GAMs under a mixed model framework. For all four methods we chose the number of interior knots as $K = 5\lceil n^{0.18}\rceil$ to satisfy Condition (C6'), where $\lceil \cdot \rceil$ denotes the ceiling function.

We simulated datasets by adapting Example 7 from the \texttt{gamSim} function available in \texttt{mgcv}, which is based on a GAM with four smoothing terms. For $i = 1,\ldots,n$, a vector of four smoothing covariates $\bm{u}_i$ was generated by simulating the $u_{i1}$ and $u_{i3}$ independently from a uniform distribution $\text{U}[0,1]$, and then generating the second and fourth elements as $u_{i2} = 0.7u_{i1} + e_1$ and $u_{i4} = 0.9u_{i3} + e_2$, where $e_1 \sim \text{U}[0,0.3]$ and $e_2 \sim \text{U}[0,0.1]$. The four smoothing functions were designed as $s_1(u) = 2\sin(\pi u), s_2(u) = \exp(2 u), s_3(u) = 0.2u^{11}\{10(1 - u)\}^6 + 10(10u)^3(1 - u)^{10}$ and $s_4(u) = 0$ \citep[see the first three panels in Figure~4.12,][]{wood2006generalized}. All smoothing terms were centered before simulating the responses. Next, we included a vector of two parametric covariates $\bm{x}_i$, where $x_{i1} = 1$ for an overall intercept and $x_{i2}$ is a binary indicator variable representing a treatment effect, say, with $x_{i2} = 1$ for the first $n/2$ observations and $x_{i2} = 0$ for the remaining $n/2$ observations. The vector of parametric coefficients was set at $\bm{\kappa} = (-1,0.5)^\top$. Finally, three types of responses generated conditional on the formulated linear predictor $\eta_i$: normal responses assuming an error variance of one, Poisson responses with a log link, and Bernoulli responses with a logit link. We considered sample sizes $n = \{100,200,500,1000\}$ and simulated 1000 datasets for each sample size  and response type. 
 
In the above simulation design, because the first three smoothing terms are not precisely a linear combination of P-splines, then it allows us to empirically assess the performance of the VA framework when the true model is not the same as (but can be well approximated by) the assumed model, thus complementing the theoretical results developed in Section~\ref{sec:asymptotic}.

We used a variety of criteria to assess performance. For overall (parametric plus smoothing) performance, we calculated the mean squared error on both the linear predictor and the mean response scale. For the parametric component with $\kappa_2 = 0.5$, we considered empirical bias and mean squared error, as well as the coverage probability and the width of the corresponding 95\% Wald confidence interval constructed using the variational information matrix in Section~\ref{subsec:parametricbit}. Finally, we assessed out of sample performance for each method by randomly removing ten observations from the dataset for validation, fitting the GAM to the remaining $(n-10)$ training observations, and constructing predictions and associated 95\% prediction intervals for the validation observations (using the methods discussed in Section~\ref{sec:inference}). Based on this, we computed the average interval width and interval score \citep{gneiting2007strictly} across the ten validation points. The interval score is calculated as follows: for a prediction interval with $100(1-\alpha)\%$ nominal coverage probability, let the lower and upper bounds of the interval be given by 
$l$ and $u$ respectively. At the ten evaluation points, the interval score is then defined as $\sum_{i=1}^{10}\left[(u_i-l_i) + 2\alpha^{-1} \mathds{1}\{y_i>u_i\} + 2\alpha^{-1} \mathds{1}\{y_i<l_i\}\right]$, where $y_i$ denotes a validation point and $\mathds{1}{\{\cdot\}}$ is the binary indicator function. 
A lower interval score corresponds to better predictive performance of the constructed prediction intervals overall i.e., the interval covers the validation point and its width is comparably small. For all assessment criteria except coverage probability, we considered both averaging and taking the median across simulated datasets, e.g., average and median mean squared errors on the linear predictor scale. Finally, the computation time in seconds for each method was also recorded. 

For brevity, we only present the main results for Poisson and Bernoulli GAMs, based on averaging across simulated datasets. The normal response GAM (which presents similar trends to the other two responses types) as well as the full simulation results (including an additional method where we considered VA assuming a block diagonal structure for $\bm{A}$ and the additional assessment criteria) are presented in the Supplementary Material. 

\subsection{Poisson Responses} \label{subsec:poissonsims}

The VA approach performed strongly compared to \texttt{mgcv} and \texttt{gamm4} (Table~\ref{tab:poisson}), with the lowest mean squared errors for the overall fit on the linear predictor scale, and the coverage probability for the treatment effect approaching the nominated 95\% level with increasing sample size. Out of sample performance, as based on lower interval scores, was also competitive for the VA approach, especially given the VA approach also produced narrower prediction intervals compared to the other procedures (see the Supplementary Material). There was little difference in the mean squared errors for the parametric component across the methods. Both implementations of \texttt{mgcv} performed somewhat poorly at the two smaller sizes, and this appears to have to be due to some instability in these algorithms for a proportion of the datasets. This, in turn, is perhaps a reflection of the instability of the penalized likelihood approach at small sample sizes, as reviewed in Section~\ref{sec:intro}. If we consider taking the median across simulated datasets, then this instability largely disappears, although the overall trends in performance were similar to those based on averaging across datasets (see the Supplementary Material).

\begin{table}[tb]
\tabcolsep 0.24in
\caption{Simulation results for Poisson GAMs, based on averaging across datasets (with sample standard deviation in parentheses for all but CI coverage). Below, we present results for the: mean squared error of the parametric component (MSE$_{\text{p}}$), coverage probability of the 95\% confidence interval for the parametric component (CI coverage$_{\text{p}}$), mean squared error for the overall fit on the linear predictor scale (MSE), and interval score for the ten out of sample validation points (Interval score).} 
\medskip \centering \label{tab:poisson}
\scalebox{0.8}{
\begin{tabular}{@{}llrrrr@{}}
    \toprule[1.5pt]
    $n$ & & VA (Unstruc) & mgcv (Default) & mgcv (P-splines) & gamm4 \\
    \cmidrule{3-6}    
      \multirow{5}{*}{100} & MSE$_{\text{p}}$ & 0.012 (0.020) & 0.010 (0.020) & 0.016 (0.026) & 0.012 (0.020) \\ 
      & CI coverage$_{\text{p}}$ & 0.972 & 0.936 & 0.957 & 0.975 \\ 
      & MSE & 0.328 (0.177) & 17.581 (186.679) & 1346.472 ($>10^4$) & 0.627 (0.515) \\ 
      & Interval score & 33.331 (5.167) & 35.469 (14.334) & 119.501 (2218.781) & 32.051 (5.365) \\ \\
      
      \multirow{5}{*}{200} & MSE$_{\text{p}}$ & 0.003 (0.004) & 0.002 (0.003) & 0.003 (0.004) & 0.003 (0.004) \\ 
      & CI coverage$_{\text{p}}$ & 0.967 & 0.941 & 0.948 & 0.961 \\ 
      & MSE & 0.205 (0.117) & 10.197 (131.503) & 36.133 (483.229) & 0.395 (0.357) \\ 
      & Interval score & 34.887 (4.776) & 35.692 (6.809) & 35.190 (8.942) & 34.169 (5.073) \\ \\

      \multirow{5}{*}{500} & MSE$_{\text{p}}$ & 0.001 (0.001) & 0.001 (0.001) & 0.001 (0.001) & 0.001 (0.001) \\ 
      & CI coverage$_{\text{p}}$ & 0.946 & 0.930 & 0.938 & 0.941 \\ 
      & MSE & 0.125 (0.085) & 2.720 (38.518) & 1145.791 $(>10^4$)& 0.206 (0.219) \\ 
      & Interval score & 36.594 (3.810) & 37.366 (3.761) & 40.420 (86.984) & 36.258 (4.027) \\ \\

      \multirow{5}{*}{1000} & MSE$_{\text{p}}$ & 0.000 (0.000) & 0.000 (0.000) & 0.000 (0.000) & 0.000 (0.000) \\ 
      & CI coverage$_{\text{p}}$ & 0.968 & 0.951 & 0.965 & 0.965 \\ 
      & MSE & 0.076 (0.053) & 0.430 (0.782) & 0.201 (1.847) & 0.111 (0.110) \\ 
      & Interval score & 37.386 (3.457) & 37.919 (3.329) & 37.310 (3.489) & 37.202 (3.651) \\ 
    \bottomrule[1.5pt]
  \end{tabular}
  }
  \end{table}

\begin{figure}[tb]
\caption{Comparative boxplots of computation time in seconds for various methods of estimating GAMs with Poisson responses. Note time on the $y$-axis is on the log scale.} \medskip \centering \label{fig:poisgam} 
\includegraphics[width = \textwidth]{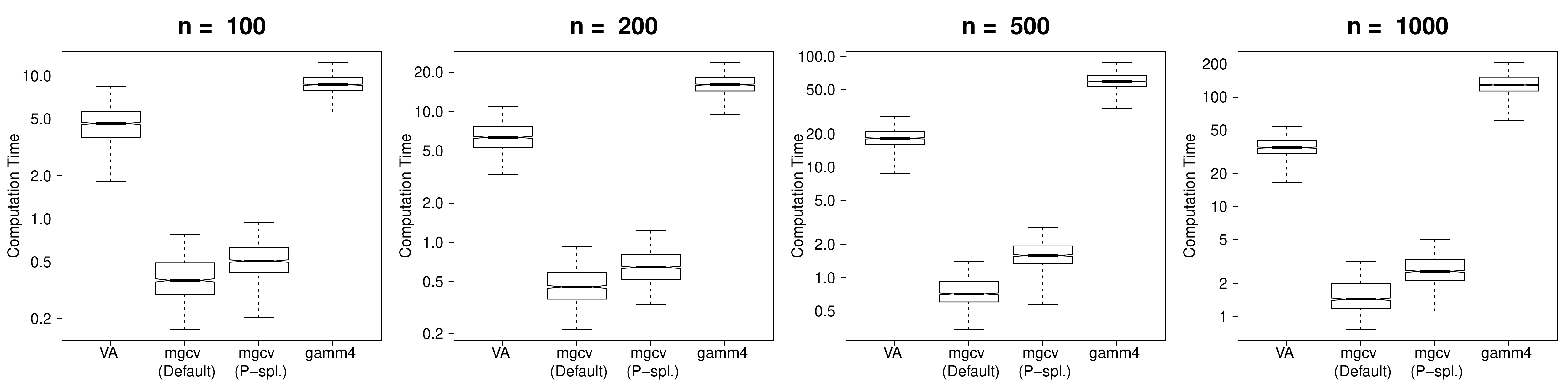}
\end{figure}

Regarding computation time, the two implementations of \texttt{mgcv} were the fastest of all the methods (Figure~\ref{fig:poisgam}). This was not surprising given that the penalized likelihood for a GAM does not involve an integral and the bulk of the computation is performed in C$++$ (our proposed VA approach was entirely implemented in \texttt{R}). Of more relevance here are the two approaches which parameterize the GAM as mixed models. Specifically, VA was both substantially faster and scaled better with sample size compared to \texttt{gamm4} (noting that the bulk of the computation for \texttt{gamm4} is also performed in C$++$.).  

\subsection{Bernoulli Responses} \label{subsec:binarysims}

The VA approach again performed competitively (Table~\ref{tab:binary}), with mean squared errors for both the overall fit and parametric component on the linear predictor scale lower than the implementations of \texttt{mgcv} and \texttt{gamm4}, and the coverage probability for the treatment effect approaching the nominated 95\% level with increasing sample size. The interval scores for VA tended to be marginally higher than \texttt{gamm4} across all sample sizes considered. Looking at the full results in the Supplementary Material, there is a suggestion that perhaps the prediction intervals were slightly too narrow, leading to slight under-coverage. Both implementations of \texttt{mgcv}, and to a lesser extent \texttt{gamm4}, performed poorly at the smaller sizes, with the instability here worse than for the Poisson responses case. This was not surprising given that Bernoulli responses contain inherently less information about the response-covariate relationships compared to Poisson responses. If we consider taking the median across simulated datasets, then this instability is still present although to a lesser extent, with the overall trends in performance similar to those based on averaging across datasets (see the Supplementary Material).
\FloatBarrier

\begin{table}[tb]
\tabcolsep 0.21in
\caption{Results for Bernoulli GAMs, based on averaging across datasets (with sample standard deviation in parentheses for all but CI coverage). Below, we present results for the: mean squared error of the parametric component (MSE$_{\text{p}}$), coverage probability of the 95\% confidence interval for the parametric component (CI coverage$_{\text{p}}$), mean squared error for the overall fit on the linear predictor scale (MSE), and interval score for the ten out of sample validation points (Interval score).} 
\medskip \centering \label{tab:binary}
\scalebox{0.8}{
\begin{tabular}{@{}llrrrrr@{}}
    \toprule[1.5pt]
    $n$ & & VA (Unstruc) & mgcv (Default) & mgcv (P-splines) & gamm4 \\
    \cmidrule{3-6}
    \multirow{5}{*}{100} & MSE$_{\text{p}}$ & 0.319 (0.492) & 2069.048 ($>10^4$) & 1468.879 (6792.697) & 5.934 (59.111) \\ 
    & CI coverage$_{\text{p}}$ & 0.964 & 0.960 & 0.990 & 0.934 \\ 
    & MSE & 2.021 (0.711) & $>10^4$ ($>10^4$) & $>10^4$ ($>10^4$) & 33.948 (306.574) \\ 
    & Interval score & 31.466 (5.696) & $>10^4$($>10^4$) & $>10^4$ ($>10^4$) & 30.305 (23.567) \\\\ 

    \multirow{5}{*}{200} & MSE$_{\text{p}}$ & 0.174 (0.263) & 712.480 ($>10^4$) & 908.577 (6071.724) & 0.269 (0.449) \\ 
    & CI coverage$_{\text{p}}$ & 0.963 & 0.928 & 0.951 & 0.931 \\ 
    & MSE & 1.191 (0.450) & $>10^4$ ($>10^4$) & $>10^4$ ($>10^4$) & 1.417 (1.150) \\ 
    & Interval score & 33.138 (5.435) & 3918.445 ($>10^4$) & $>10^4$ ($>10^4$) & 29.235 (6.424) \\\\
    
    \multirow{5}{*}{500} & MSE$_{\text{p}}$ & 0.067 (0.093) & 0.083 (0.114) & 0.082 (0.115) & 0.083 (0.115) \\ 
    & CI coverage$_{\text{p}}$ & 0.962 & 0.956 & 0.958 & 0.952 \\ 
    & MSE & 0.635 (0.244) & 17.013 (370.524) & $>10^4$ ($>10^4$) & 0.650 (0.385) \\ 
    & Interval score & 35.086 (4.581) & 33.735 (5.367) & 37.202 (80.163) & 32.158 (5.704) \\\\ 

    \multirow{5}{*}{1000} & MSE$_{\text{p}}$ & 0.037 (0.053) & 0.041 (0.059) & 0.039 (0.057) & 0.042 (0.061) \\ 
    & CI coverage$_{\text{p}}$ & 0.952 & 0.961 & 0.957 & 0.952 \\ 
    & MSE & 0.356 (0.142) & 0.619 (2.702) & 0.566 (3.509) & 0.350 (0.176) \\ 
    & Interval score & 35.932 (4.385) & 34.997 (4.801) & 33.685 (5.352) & 33.765 (5.268) \\     
    \bottomrule[1.5pt]
  \end{tabular}
  }
  \end{table}

\begin{figure}[tb]
\caption{Comparative boxplots of computation time in seconds for various methods of estimating GAMs with Bernoulli responses. Note time on the $y$-axis is on the log scale.} \medskip \centering \label{fig:binarygam}
\includegraphics[width = \textwidth]{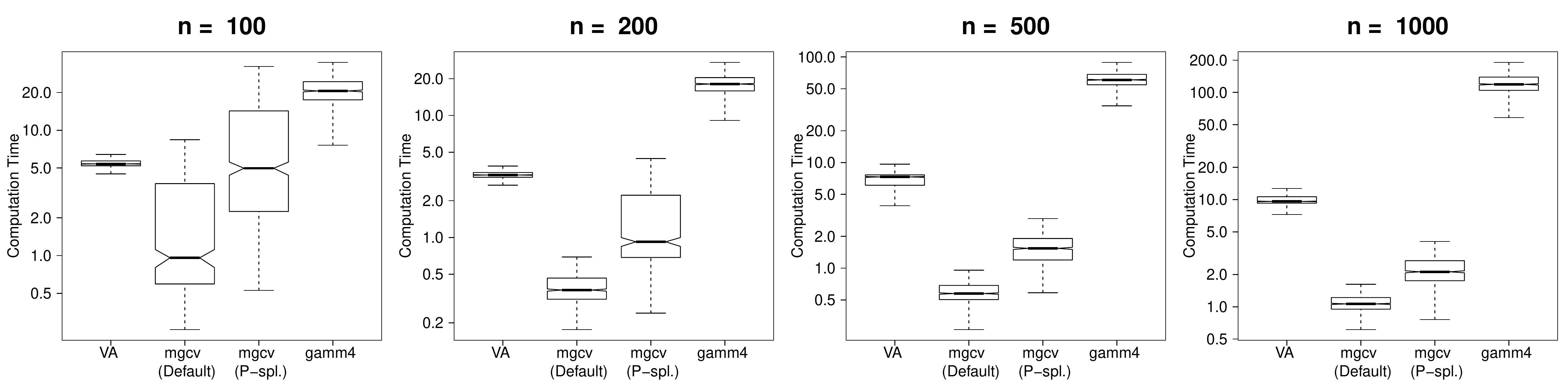}
\end{figure}

Results for computation time were were similar to those for Poisson GAMs in Section~\ref{subsec:poissonsims} (Figure~\ref{fig:binarygam}). In particular, of the two mixed model approaches to GAMs, the VA approach was substantially faster and scaled better with sample size compared to \texttt{gamm4}. 

In the full simulation results presented in the Supplementary Material, we found that the VA approach assuming a suboptimal block diagonal structure for $\bm{A}$ had little impact on performance, but turned out to be slower than assuming an unstructured form. We suspect that this is due to the way the updates for $\bm{A}$ were coded, i.e., under the block diagonal structure, we used a for loop was used to update each submatrix $\bm{A}_j$, whereas the unstructured involved a single update of the full $\bm{A}$. The Supplementary Material also presents results for the normal response case which, not surprisingly, showed little difference in performance between all the methods considered. Interestingly, even in cases where the marginal likelihood is tractable, the VA approach was still faster and scaled considerably better than \texttt{gamm4}.


\section{Application to Union Membership Dataset} \label{sec:application} 

We applied the VA approach for fitting GAMs to a 1985 North American population survey containing information on union membership and various worker's attributes, with the aim being to uncover associations between workers' characteristics and their probability of union membership \citep{berndt1991practice,ruppert2003semiparametric}. The data consist of $n=534$ observations, with the response being a Bernoulli variable of whether they were a member of a union (1 $=$ yes; 0 $=$ no), and six covariates: gender (1 $=$ female; 0 $=$ male), race (1 $=$ white; 0 $=$ other), an indicator variable for whether the worker lives in the south (1 $=$ yes; 0 $=$ no), age in years, hourly wage, and number of years in education. 

\begin{figure}[t]
\caption{Smooths from the fitted GAM using the VA approach, regressing union membership as a function of six covariates. Results show no evidence of a relationship between probability of a worker being in a union and their education duration or age, but a strong non-linear relationship with their hourly wage.} \label{fig:wagevagam} \centering \medskip
\includegraphics[width = \textwidth]{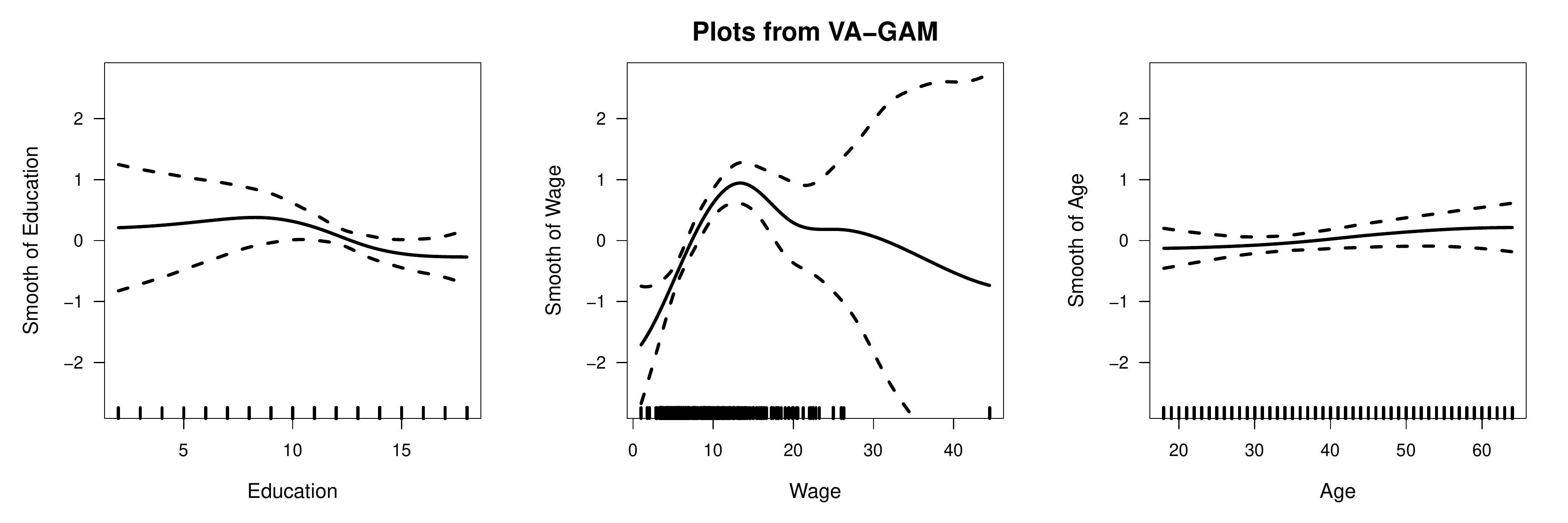}
\end{figure}

We fitted a GAM to union membership, assuming a Bernoulli response with logit link, with the first three binary covariates (gender, race, south) entered as parametric terms and the last three continuous covariates (age, wage, education) entered as smooth terms. Following some exploratory analysis, we found that $K = 8$ interior knots were sufficient for all three continuous covariates. Indeed, when we fitted a GAM to this dataset using \texttt{mgcv} it also selected $K = 8$, and so for consistency we chose to use the same number of interior knots. Additional exploration using additional knots showed little difference in the fitted smooths. Inference on the parametric component of the fitted GAM showed evidence that females workers ($\hat{\kappa} = -0.700$; 95\% Wald confidence interval: $[-1.216,-0.186]$) and workers of white descent ($\hat{\kappa} = -0.724$; 95\% Wald confidence interval: $[-1.306,-0.142]$) were less likely to be a union member. However, at the 5\% significance level there was no strong evidence of an association between union membership and living in the south ($\hat{\kappa} = -0.498$; 95\% Wald confidence interval: $[-1.074,0.079]$). For the three continuous covariates, plots of the resulting smooths and associated pointwise confidence bands suggested there was only a strong, non-linear relationship between probability of union membership and hourly wage (Figure~\ref{fig:wagevagam}). This result was confirmed in the Wald tests for the smooth components, where there was no evidence of an association to age or education ($p$-values for both exceeding 0.4), but strong evidence of an association to wage ($p$-value $< 0.01$). The resulting smooth suggested that the probability union membership increased until the hourly wage hit around $\$15$, then decreased steeply until around $\$22$. 

The conclusions obtained above from using the VA framework were very similar to those obtained when GAMs with the same setup were fitted using the \texttt{mgcv} and \texttt{gamm4} packages (see Supplementary Material), with the only noteworthy difference being that the latter two fitted models revealed borderline evidence of the smooth curve for education being significant ($p$-values for both \texttt{mgcv} and \texttt{gamm4} were around 0.05). On computation time, \texttt{mgcv} was the fastest of the methods (0.230 seconds), although the VA method (1.838 seconds) was faster than \texttt{gamm4} (5.599 seconds). 

\section{Discussion} \label{sec:discuss}

In this article, we have proposed a variational approximations framework for GAMs, which offers a balance between speed and stability: VA is closer to the penalized likelihood approach in terms of computational speed, but it allows us to borrow the inference tools and generally stronger finite sample performance of the mixed model framework. We derived fully closed-form approximations to the marginal likelihood for commonly observed responses and demonstrated how they could be maximized straightforwardly. We proposed how several key aspects of inference can be performed for GAMs using the VA approach, including a variational information matrix for inference on the parametric component and fully closed-form updates of the smoothing parameters. Simulation studies showed the VA approach performed similarly to or better than current approaches for estimation and inference with GAMs, while achieving speeds competitive to that of the penalized likelihood approach.

It is important to stress that the consistency and asymptotic normality we develop here for the VA estimates are driven not so much by the fact that the approximation to the marginal log-likelihood is asymptotically negligible, but because, as discussed in Section~\ref{sec:asymptotic}, we work under the framework that the true model is a GLM with the relevant spline basis functions, and that the effect of the penalization due to smoothing is asymptotically dominated by the likelihood of the data. Studying the error produced by applying variational approximations to the marginal likelihood and other related quantities such as the information matrix as in (\ref{eqn:infomat}) remains an open and challenging task of interest \citep{blei2017variational}. On the other hand, the fact that we can derive appealing asymptotic results without explicitly studying the approximation error is an important advance and promotes a strong theoretical argument for using VA in semiparametric regression; we refer the reader to \citet{westling2017} and \citet{chen2018bootstrap} for some examples of recent theoretical developments in large sample theory variational approximations in statistics, neither of which involve studying the approximation error directly.

This article makes an important first step into using VA for semiparametric regression, and opens up a myriad of research opportunities. These include but are not limited into VA for generalized additive mixed models, other more sophisticated basis functions that can accommodate interactions more straightforwardly such as thin plate regression splines, adaptive quadratic penalties with more than one smoothing parameter, settings where both the number of covariates to be smoothed $q$ as well as basis dimension per covariate $d$ grows where penalties may be required to choose the basis functions (see \citet{xue09} and also \citet{ormerod2016variational} for an example of related work in the VB setting), and other types of semiparametric regression models such as single index models. With the ever-increasing volume of data being collected, the VA framework seems to provide a good compromise between stability and speed for estimation and inference in semiparametric regression.

\section*{Acknowledgements}
Thanks to Wade Blanchard and Alan Welsh for useful discussions. SM and FKCH were both partly supported by Australian Research Council discovery project grant DP180100836. HLS was partially supported by a research school grant from the Australian National University.

\bibliographystyle{apalike} 
\bibliography{thesis}

\begin{thebibliography}{}

\bibitem[Berndt, 1991]{berndt1991practice}
Berndt, E. (1991).
\newblock {\em The Practice of Econometrics: Classic and Contemporary}.
\newblock Addison-Wesley Publishing Company, Reading, Massachusetts.

\bibitem[Blei et~al., 2017]{blei2017variational}
Blei, D.~M., Kucukelbir, A., and McAuliffe, J.~D. (2017).
\newblock {Variational inference: A review for statisticians}.
\newblock {\em Journal of the American Statistical Association}, In press.

\bibitem[Blei and Lafferty, 2007]{blei07}
Blei, D.~M. and Lafferty, J.~D. (2007).
\newblock A correlated topic model of science.
\newblock {\em The Annals of Applied Statistics}, 1:17--35.

\bibitem[Breslow and Clayton, 1993]{breslow1993approximate}
Breslow, N.~E. and Clayton, D.~G. (1993).
\newblock {Approximate inference in generalized linear mixed models}.
\newblock {\em Journal of the American Statistical Association}, 88:9--25.

\bibitem[Chen et~al., 2018]{chen2018bootstrap}
Chen, Y.~C., Wang, Y.~S., and Erosheva, E.~A. (2018).
\newblock {On the use of bootstrap with variational inference: Theory,
  interpretation, and a two-sample test example}.
\newblock {\em The Annals of Applied Statistics}, In press.

\bibitem[Dempster et~al., 1977]{dempster77}
Dempster, A.~P., Laird, N.~M., and Rubin, D.~B. (1977).
\newblock {Maximum likelihood from incomplete data via the EM algorithm}.
\newblock {\em Journal of the Royal Statistical Society: Series B (Statistical
  Methodology)}, 39:1--38.

\bibitem[Eilers and Marx, 1996]{eilers1996}
Eilers, P. H.~C. and Marx, B.~D. (1996).
\newblock {Flexible smoothing with $B$-splines and penalties}.
\newblock {\em Statistical Science}, 11:89--121.

\bibitem[Fan and Peng, 2004]{fan04}
Fan, J. and Peng, H. (2004).
\newblock {Nonconcave penalized likelihood with a diverging number of
  parameters}.
\newblock {\em The Annals of Statistics}, 32:928--961.

\bibitem[Gneiting and Raftery, 2007]{gneiting2007strictly}
Gneiting, T. and Raftery, A.~E. (2007).
\newblock {Strictly proper scoring rules, prediction, and estimation}.
\newblock {\em Journal of the American Statistical Association}, 102:359--378.

\bibitem[Hall et~al., 2011]{hall2011asymptotic}
Hall, P., Pham, T., Wand, M.~P., and Wang, S.~S. (2011).
\newblock {Asymptotic normality and valid inference for Gaussian variational
  approximation}.
\newblock {\em The Annals of Statistics}, 39:2502--2532.

\bibitem[Hastie and Tibshirani, 1986]{hastie86}
Hastie, T. and Tibshirani, R. (1986).
\newblock {Generalized additive models}.
\newblock {\em Statistical Science}, 1:297--310.

\bibitem[Hui et~al., 2017a]{hui16crepe}
Hui, F. K.~C., Mueller, S., and Welsh, A.~H. (2017a).
\newblock Hierarchical selection of fixed and random effects in generalized
  linear mixed models.
\newblock {\em Statistica Sinica}, 27:501--518.

\bibitem[Hui et~al., 2017b]{hui2016variational}
Hui, F. K.~C., Warton, D.~I., Ormerod, J.~T., Haapaniemi, V., and Taskinen, S.
  (2017b).
\newblock Variational approximations for generalized linear latent variable
  models.
\newblock {\em Journal of Computational and Graphical Statistics}, 26:35--43.

\bibitem[Kauermann, 2005]{kauermann2005note}
Kauermann, G. (2005).
\newblock A note on smoothing parameter selection for penalized spline
  smoothing.
\newblock {\em Journal of statistical planning and inference}, 127:53--69.

\bibitem[Kauermann et~al., 2009]{kauermann2009some}
Kauermann, G., Krivobokova, T., and Fahrmeir, L. (2009).
\newblock {Some asymptotic results on generalized penalized spline smoothing}.
\newblock {\em Journal of the Royal Statistical Society: Series B (Statistical
  Methodology)}, 71:487--503.

\bibitem[Knowles and Minka, 2011]{knowles11}
Knowles, D.~A. and Minka, T. (2011).
\newblock {Non-conjugate variational message passing for multinomial and binary
  regression}.
\newblock In Shawe-Taylor, J., Zemel, R.~S., Bartlett, P.~L., Pereira, F., and
  Weinberger, K.~Q., editors, {\em Advances in Neural Information Processing
  Systems 24}, pages 1701--1709. Curran Associates, Inc.

\bibitem[Krivobokova et~al., 2010]{krivobokova2010}
Krivobokova, T., Kneib, T., and Claeskens, G. (2010).
\newblock Simultaneous confidence bands for penalized spline estimators.
\newblock {\em Journal of the American Statistical Association}, 105:852--863.

\bibitem[Lin and Zhang, 1999]{lin1999inference}
Lin, X. and Zhang, D. (1999).
\newblock {Inference in generalized additive mixed models by using smoothing
  splines}.
\newblock {\em Journal of the Royal Statistical Society. Series B (Statistical
  Methodology)}, 61:381--400.

\bibitem[Louis, 1982]{louis82}
Louis, T.~A. (1982).
\newblock {Finding the observed information matrix when using the EM
  algorithm}.
\newblock {\em Journal of the Royal Statistical Society: Series B (Statistical
  Methodology)}, 44:226--233.

\bibitem[Lu et~al., 2008]{lu2008penalized}
Lu, Y., Zhang, R., and Zhu, L. (2008).
\newblock Penalized spline estimation for varying-coefficient models.
\newblock {\em Communications in Statistics--Theory and Methods},
  37:2249--2261.

\bibitem[Luts et~al., 2014]{luts2014real}
Luts, J., Broderick, T., and Wand, M.~P. (2014).
\newblock Real-time semiparametric regression.
\newblock {\em Journal of Computational and Graphical Statistics}, 23:589--615.

\bibitem[Luts and Wand, 2015]{luts2015variational}
Luts, J. and Wand, M.~P. (2015).
\newblock Variational inference for count response semiparametric regression.
\newblock {\em Bayesian Analysis}, 10:991--1023.

\bibitem[Ormerod and Wand, 2010]{ormerod2010explaining}
Ormerod, J. and Wand, M. (2010).
\newblock {Explaining variational approximations}.
\newblock {\em The American Statistician}, 64(2):140--153.

\bibitem[Ormerod and Wand, 2012]{ormerod2012gaussian}
Ormerod, J. and Wand, M. (2012).
\newblock {Gaussian variational approximate inference for generalized linear
  mixed models}.
\newblock {\em Journal of Computational and Graphical Statistics}, 21:2--17.

\bibitem[Ormerod et~al., 2016]{ormerod2016variational}
Ormerod, J.~T., You, C., and Mueller, S. (2016).
\newblock {A variational Bayes approach to variable selection}.
\newblock {\em Electronic Journal of Statistics}, In review.

\bibitem[Reiss and Ogden, 2009]{reiss2009smoothing}
Reiss, P.~T. and Ogden, R.~T. (2009).
\newblock Smoothing parameter selection for a class of semiparametric linear
  models.
\newblock {\em Journal of the Royal Statistical Society: Series B (Statistical
  Methodology)}, 71:505--523.

\bibitem[Rigby and Stasinopoulos, 2005]{rigby05}
Rigby, R.~A. and Stasinopoulos, D.~M. (2005).
\newblock Generalized additive models for location, scale and shape,(with
  discussion).
\newblock {\em Journal of the Royal Statistical Society. Series C (Applied
  Statistics)}, 54:507--554.

\bibitem[Ruppert, 2002]{ruppert2002}
Ruppert, D. (2002).
\newblock {Selecting the number of knots for penalized splines}.
\newblock {\em Journal of Computational and Graphical Statistics}, 11:735--757.

\bibitem[Ruppert et~al., 2003]{ruppert2003semiparametric}
Ruppert, D., Wand, M., and Carroll, R. (2003).
\newblock {\em Semiparametric Regression}.
\newblock Cambridge University Press, New York.

\bibitem[Wang et~al., 2011]{wang2011asymptotics}
Wang, X., Shen, J., Ruppert, D., et~al. (2011).
\newblock On the asymptotics of penalized spline smoothing.
\newblock {\em Electronic Journal of Statistics}, 5:1--17.

\bibitem[Westling and McCormick, 2017]{westling2017}
Westling, T. and McCormick, T.~H. (2017).
\newblock {Consistency, calibration, and efficiency of variational inference}.
\newblock {\em arXiv preprint}, 1510.08151v3.

\bibitem[Wood, 2006]{wood2006generalized}
Wood, S. (2006).
\newblock {\em {Generalized Additive Models: An Introduction with R}}.
\newblock CRC press, Boca Raton, FL.

\bibitem[Wood and Scheipl, 2016]{gamm4}
Wood, S. and Scheipl, F. (2016).
\newblock {\em gamm4: Generalized Additive Mixed Models using `mgcv' and
  `lme4'}.
\newblock R package version 0.2-4. URL:
  \url{https://CRAN.R-project.org/package=gamm4}.

\bibitem[Wood, 2011]{wood2011fast}
Wood, S.~N. (2011).
\newblock Fast stable restricted maximum likelihood and marginal likelihood
  estimation of semiparametric generalized linear models.
\newblock {\em Journal of the Royal Statistical Society: Series B (Statistical
  Methodology)}, 73:3--36.

\bibitem[Wood et~al., 2016]{wood2016}
Wood, S.~N., Pya, N., and S\"{a}fken, B. (2016).
\newblock Smoothing parameter and model selection for general smooth models.
\newblock {\em Journal of the American Statistical Association},
  111:1548--1563.

\bibitem[Xue, 2009]{xue09}
Xue, L. (2009).
\newblock {Consistent variable selection in additive models}.
\newblock {\em Statistica Sinica}, 19:1281--1296.

\bibitem[Yu and Ruppert, 2002]{yu2002penalized}
Yu, Y. and Ruppert, D. (2002).
\newblock Penalized spline estimation for partially linear single-index models.
\newblock {\em Journal of the American Statistical Association}, 97:1042--1054.

\end{thebibliography}

\end{document}